\documentclass[11pt,
]{article}                          

\usepackage{amssymb,amsfonts,amstext,amsmath,amsthm,latexsym,mathrsfs}
\usepackage{mathbbol}

\usepackage{mparhack}
\usepackage{makeidx}

\makeindex

\makeatletter
\if@twoside%
\else\fi%
\makeatother  

\input{coh.sty}

\begin{document}

\title
{In Cohen generic extension, 
every countable OD set of reals belongs to the ground model
}

\author 
{
Vladimir~Kanovei\thanks{IITP RAS and MIIT,
  Moscow, Russia, \ {\tt kanovei@googlemail.com} --- contact author. 
}  
}

\date 
{\today}

\maketitle

\begin{abstract}
It is true in the Cohen generic extension of $\rL$, 
the constructible universe, that every countable 
ordinal-definable set of reals belongs to $\rL$.
\end{abstract}

\bte
\lam{cohT}
Let\/ $a\in\bn$ be a Cohen-generic real over\/ $\rL$.
Then it is true in\/ $\rL[a]$ that if\/ $X\sq\bn$ is 
a countable\/ $\od$ set then\/ $X\in\rL$.
\ete

One may expect such a result of any homogeneous forcing notion. 
For instance, Theorem~\ref{cohT} is true for the Solovay model 
(the extension of $\rL$ by Levy cardinal collapse up to 
an inaccessible cardinal \cite{solo})  
--- but by a different argument. 
One hardly can doubt that any typical homogeneous extension 
(Solovay-random, Sacks, Hehler, and the like) also satisfies 
the same result, but it's not easy to manufacture a proof of 
sufficient generality.

On the contrary, non-homogeneous forcing notions may lead 
to models with countable $\od$ non-empty sets of reals with 
no $\od$ elements \cite{kl:cds}, and such a set can even have 
the form of   
a $\ip12$ \dd\Eo equivalence class \cite{kl:dec}.

\bpf
Let $C=\nse$ be the Cohen forcing. 
First of all, it suffices to prove that 
(it is true in\/ $\rL[a]$ that) if\/ $X\sq\bn$ is 
a countable\/ $\od$ set then\/ $X\sq\rL$.
Indeed, as the Cohen forcing is homogeneous, any 
statement about sets in $\rL$, the ground model, 
is decided by the weakest condition.

There is a formula $\vpi(x)$ 
with an unspecified ordinal $\al_0$ as a parameter, 
such that $X=\ens{x\in\bn}{\vpi(x)}$ in $\rL[a]$, and then 
there is a condition $p_0\in C$ such that $p_0\su a$ and 
$p_0$ \dd Cforces that 
$\ens{x\in\bn}{\vpi(x)}$ is a countable set.
Suppose to the contrary that $X\not\sq\rL$, so that $p_0$ 
also forces $\sus x\:(x\nin\rL\land\vpi(x))$.

There is a sequence $\sis{t_n}{n<\om}\in\rL$ of \dd Cnames, 
such that if $b\in\bn$ is Cohen generic and $p_0\su b$ then 
it is true in $\rL[b]$ that 
$\ens{x\in\bn}{\vpi(x)}=\ens{\bint{t_n}b}{n<\om}$, 
where $\bint{t}x$ is the interpretation of a \dd Cname $t$ 
by a real $x\in\bn.$
Let $T\in\rL$ be the \dd Cname for 
$\ens{\bint{t_n}{\dot a}}{n<\om}$.
Thus we assume that $p_0$ forces 
$$
\bint{T}{\dot a}= \ens{\bint{t_n}{\dot a}}{n<\om}=
\ens{x\in\bn}{\vpi(x)}
\not\sq\rL\,,
\eqno(1)
$$
where $\dot a$ is the canonical name for the \dd Cgeneric real.

Let  
$\dal,\dar$ be canonical \dd{(C\ti C)}names for the left, 
resp., right of the terms of a\/ \dd{(C\ti C)}generic pair 
of reals\/ $\ang{\alev,\apra}$.

\bcor
\lam{cor}
The pair\/ $\ang{p_0,p_0}$ \dd{(C\ti C)}forces 
over\/ $\rL$ that\/ 
$\bint{T}{\dot a\lef}\ne \bint{T}{\dot a\rig}$.
\ecor 
\bpf
$\rL[\alev]\cap\rL[\apra]\cap\bn\sq\rL$ due to the mutual 
genericity of $\alev,\apra$.
\epf 

Now pick a regular cardinal $\ka>\al_0$. 
\vyk{
By a simple forcing argument, it follows from the above 
that if $\ang{\alev,\apra}$ is a \dd{(C\ti C)}generic pair 
over $\rL_\ka$ with $p_0\su \alev$, $p_0\su\apra$, then still 
$\bint{T}{\dot a\lef}\ne \bint{T}{\dot a\rig}$. 
}%
Consider, in $\rL$, a countable submodel $\gM$ of $\rL_\ka$ 
containing $\al_0$ and all names $t_n$ and $T$. 
Let $\pi:\gM\to\bgm$ be the Mostowski collapse onto a 
transitive set $\bgm$. 

\bcor
\lam{cor'}
It is true in\/ 
$\bgm$ 
that\/ 
$\ang{p_0,p_0}$ 
\dd{(C\ti C)}forces\/ 
$\bint{T}{\dot a\lef}\ne \bint{T}{\dot a\rig}$.
\ecor 
\bpf
By the elementarity, this holds in $\gM$. 
Further we have $\pi(t_n)=t_n$ and $\pi(T)=T$ because 
the names $t_n$ and $T$ belong to the transitive part 
of $\gM$.
\epf 

\bcor
\lam{cor"}
If\/ $\ang{\alev,\apra}$ 
is a\/ \dd{(C\ti C)}generic pair 
over\/ $\bgm$ with $p_0\su \alev$, $p_0\su\apra$, 
then\/ 
$\bint{T}{\alev}\ne \bint{T}{\apra}$.\qed 
\ecor 

By the countability, there is a real $z\in\bn\cap\rL$ 
satisfying $z(j)=0$ for all $j<\dom{p_0}$ and 
\dd Cgeneric over $\bgm$, so that $\bgm[z]$ is a set in $\rL$.
Let $x\in\bn$ be \dd Cgeneric over $\rL$, with $p_0\su x$. 
Then, as $z\in\rL$, the real $y$ defined by $y(k)=z(k)+x(k)$, 
$\kaz k$, is \dd Cgeneric over $\rL$ as well, and we have 
$\rL[x]=\rL[y]$ and still $p_0\su y$. 
It follows from (1) that 
$\bint{T}{\dot a\lef}=\bint{T}{\dot a\rig}$ 
(an \od\ set of reals in $\rL[x]=\rL[y]$). 

But on the other hand by the product forcing theorem 
and the choice of $z$, the pair $\ang{x,y}$ is 
\dd{(C\ti C)}generic over $\bgm$, and hence 
$\bint{T}{\dot a\lef}\ne \bint{T}{\dot a\rig}$ 
by Corollary~\ref{cor"}, which is a contradiction. 
\epF{Theorem~\ref{cohT}}

\bre
\lam{r1}
The Solovay model \cite{solo} admits a somewhat stronger 
result established in \cite{arXiv1}, namely, 
any countable non-empty $\od$ 
{\it set of sets of reals\/} consists of $\od$ elements 
(sets of reals). 
We don't know whether this is true in the Cohen generic extension 
$\rL[a]$. 
\ere

\bre
\lam{r2}
Is Theorem~\ref{cohT} true for other popular 
forcing notions like \eg\ 
the random forcing? 
The proof above crucially employs the countability of 
the Cohen forcing.
\ere

\bibliographystyle{plain}
{\small

%
}

\end{document}